\documentclass{amsart}
\usepackage{amssymb,amsmath,amsthm,epsf,epsfig,dsfont,bbm}

\usepackage{latexsym}
\usepackage{upref, eucal}
\usepackage[all]{xy}

\newcommand{\ch}{\mb{char }}
\newcommand{\vpi}{v_\pi}
\newcommand{\oexp}{\overline{\exp}_\d}

\newcommand {\nc} {\newcommand}
\newcommand {\enm} {\ensuremath}

\DeclareMathOperator{\Spec}{\mathrm{Spec}}

\def \d{\delta}

\nc {\bdm} {\begin{displaymath}}
\nc {\edm} {\end{displaymath}}

\newtheorem {theorem} {\bf{Theorem}}[section]
\newtheorem {lemma}[theorem] {\bf Lemma}
\newtheorem {proposition}[theorem] {\bf Proposition}

\newtheorem {corollary}[theorem] {\bf Corollary}
\numberwithin {equation}{section}

\newcommand{\Ou}{\enm{\mathcal{O}}}

\nc{\J}{\enm{\mathcal{J} }}
\nc {\Z} {\enm{\mathbb{Z}}}
\nc {\form}[1] {\enm{\mbox{\underline{for}}}_{#1}}
\nc {\prol}[1] {\enm{\mbox{\underline{prol}}_{{#1}^*}}}

\nc {\stk} {\stackrel}

\newcommand{\map}{\rightarrow}

\newcommand{\beqar}{\begin{eqnarray*}}
\newcommand{\eeqar}{\end{eqnarray*}}
\newcommand{\inj}{\hookrightarrow}
\newcommand{\Pn}[2] {\ensuremath{ {\mathbb{P}}^{#1}_{#2}}}
\nc{\Quot}[3]{\enm{ {\mathfrak{Quot}_{ {#1}/{#2}/{#3}}}}}
\nc{\Hilb}[2]{\enm{ {\mathfrak{Hilb}_{ {#1}/{#2}}}}}
\newcommand{\mfrak}[1]{\mathfrak{#1}}

\newcommand{\bb}[1]{\mathbb{#1}}

\nc {\Coh}[4] {\ensuremath{H^{#1}(\Pn{#2}{},{#3}({#4}))}}
\nc {\Ch}[3] {\enm{H^{#1}(X_t,{#2}_t({#3}))}}
\nc {\Qphi}[4]{\enm{ {\mathfrak{Quot}^{~#4}_{ {#1}/{#2}/{#3}}}}}
\nc {\Gra}[4]{\enm{ {\mathfrak{Grass}_{#2}({#3},{#4})}}}
\nc {\HomA}[2]{\enm{\mathrm{Hom}_A{#1}{#2}}}
\nc {\tr}{\mathrm{tr}}

\nc {\C}[2]{\enm{\left(\begin{array}{l} {#1} \\ {#2} \end{array} \right)}}
\nc {\mat}[4]{\enm{\left(\begin{array}{ll}{#1} & {#2} \\ {#3} & {#4}
\end{array}\right)}}

\def \mb{\mbox}

 \def \Z{{\mathbb Z}}

   \def \h{\hat{\ }}

\def \d{\delta} \def \bZ{{\mathbb Z}}  \def \bF{{\bf F}}

 \def \bF{{\mathbb F}}

\def \hA{\hat{A}}

\def \R1{R((q))[q']\h}

\usepackage{amsmath}
\usepackage{mathtools}

\DeclareMathOperator{\texp}{\overline{\exp}_\d}

\begin{document}
\title{Lift of Frobenius and Descent to Constants}
\author{Arnab Saha}
\maketitle

\begin{abstract}
In differential algebra, a proper scheme $X$ defined over an algebraically
closed field $K$
with a derivation $\partial$ on it descends to the field of constants 
$K^\partial$ if $X$ itself lifts the derivation $\partial$. This is a result by
A. Buium. Now in the arithmetic case, the notion of a derivation is replaced
by the notion of a $\pi$-derivation $\d$ or equivalently in the flat case, 
a lift of Frobenius $\phi$. We will show an analogous result in the 
arithmetic case of equal characteristic.
We show our results using the 
arithmetic analogue of Taylor expansion using Witt vectors.
\end{abstract}

\section{Introduction}
In \cite{bui87}, Buium shows that given a proper smooth scheme $X$ over an
algebraically closed field $K$ with a
derivation $\partial$ on it, if the structure sheaf $\Ou_X$ has a derivation
lifting the one on $K$, then $X$ descends to the field of constants $K^\partial
:= \{x \in K~|~ \partial x =0\}$. Consequently, Gillet gave another 
proof of this result \cite{gil02}. In this paper, we will prove an analogous 
result in the arithmetic case where the notion of derivation is replaced
by a $\pi$-derivation $\d$ or equivalently in the flat case, by a lift of 
Frobenius $\phi$. Before we state our main results, we would like to give a 
brief background.

By a tuple of $(B,\mfrak{p})$ we will understand the following data:
$B$ is a Dedekind domain, $\mfrak{p}$ a fixed non-zero prime of $B$ with 
$k:=B/\mfrak{p}$
a finite field and let $q = |k|$. Then the identity map $\mathbbm{1}:B \map B$
is a $q$-power lift of Frobenius since for all $x \in B$, $x \equiv x^q$ 
mod $\mfrak{p}$.
Let $R$ be the $\mfrak{p}$-adic completion of $B$ and $\iota:B \map R$ be the 
canonical injective map. Also let $\pi \in B$ be such that $\iota(\pi)$ 
is a generator of the maximal ideal in $\mfrak{m} \subset R$.
Since $\iota$ is an injection, we will sometimes, by abuse of notation, 
consider $\pi$ as an element of $R$ as well. We will call $B$ or $R$ is of equal
characteristic if $\ch B = \ch R = \ch k$. Let $S = \Spec B$.

Let $A$ be a $R$-algebra. Then as in \cite{bor11a} one can 
define the $\pi$-typical Witt vectors $W(A)$ with respect to $R$.
For example, when $R=
\bZ_p$ and $\mfrak{m}=(p)$ for some prime $p$, then $W(A)$ are the usual 
$p$-typical Witt vectors.

Let $A$ be a flat seperated $\pi$-adically complete $R$-algebra with a $q$-power
Frobenius $\phi$ which is identity on $R$. Then one can consider an operator
called the $\pi$-derivation $\d$ on $A$ associated to $\phi$ given by
$\d x = \frac{\phi(x)-x^q}{\pi}$ for all $x \in X$.
By the universal 
property of Witt vectors, we obtain a canonical map $\exp_\d:A \map W(A)$ 
given by 
\begin{align}
\exp_\d(x)=(P_0(x),P_1(x),\cdots) 
\end{align}
where $P_0(x)=x$, $P_1(x)=\d x$. This map should be viewed as the analogue of 
the {\it Hasse-Schmidt} map, $\exp_\partial: A \map D(A):=A[[t]]$  in the case 
of usual derivation $\partial$ in
differential algebra, given by
\begin{align}
\exp_\partial(x)= \sum_{i=0}^\infty \frac{\partial^{(i)}x}{i!}
\end{align}
Let $A_0:=A/\pi A$ and $u$ be the quotient map $u:A \map A_0$. Then 
consider the composition $\overline{\exp}_\d:
A \map W(A_0)$ given by $\overline{\exp}_\d:= W(u) \circ \exp_\d$. We say
that $\overline{\exp}_\d$ is the {\it arithmetic Taylor expansion 
centered at $(\pi)$} and is an arithmetic analogue of the usual Taylor 
expansion principle-- this is discussed in detail in section $3$. 

Analogous to differential algebra, we
define the set of $\d$-constants of $A$ as $A^{\d}:= \{x \in A ~|~ 
\d x = 0\}$. Then it is easy to see that it is a multiplicatively closed set. 
However if $R$ is of equal characteristic,
then $A^{\d}$ is also closed under addition and 
$A^{\d}$ is then a subring of $A$. Our first result, theorem 
\ref{nokernel}, states that the arithmetic Taylor expansion 
$\overline{\exp}_\d$ is injective. 

\begin{theorem}
\label{nokernel}
Let $A$ be a separated flat $R$-algebra with a $\pi$-derivation $\d$ 
(equivalently a $q$-power Frobenius $\phi$) on it.  Then the arithmetic 
Taylor expansion map $\overline{\exp}_\d: A \map W(A_0)$ is an injection.
\end{theorem}

As a consequence we obtain corollary \ref{topology} which shows
that the $\pi$-adic topology on such an $A$ with a $\pi$-derivation $\d$ on it,
is obtained by the
restriction of the topology on $W(A_0)$ induced from the system of open
neighbourhoods $\{I_n\}_n$ around $0$, where for each $n\geq 1$, 
$$I_n:=\{x \in W(A_0)~|~ x=(x_0,x_1, \cdots ), x_i \in A_0 \mb{ and }
x_0= x_1= \cdots =x_{n-1} = 0\}$$

\begin{theorem}
\label{iso}
Assume that $R$ is of equal characteristic.
Let $A$ be a separated flat $\pi$-adically complete $R$-algebra
 and with a lift of Frobenius $\phi$ which induces an 
isomorphism on $A_0$. Then $A^\d \simeq A_0$.

Furthermore we have $A \simeq A^\d[[\pi]] \simeq A_0[[\pi]]$. 
\end{theorem}

In the above result, we would like to remark that it follows from the theory
of Witt vectors that such an $A$ is isomorphic to $W(A_0)\simeq A_0[[\pi]]$. 
The key point of the result is to show that the subring of constants $A^\d$ is 
isomorphic to $A_0$.



Assume further that $B$ is of equal characteristic and
is also a $k$-algebra where recall $k=B/\mfrak{p}$. Let $X$ be a scheme over
$S$ such that its structure sheaf $\Ou_X$ has a lift of $q$-power Frobenius 
$\phi$ which restricts to identity on $\Ou_S$. Let $\overline{X}= X \times_S
(\Spec k)$.
Furthermore, we will say $\phi$ is perfect if $\phi$ induces an automorphism 
on $\Ou_{\overline{X}}$.

Then as an application of Theorem \ref{iso} we prove the following analogous 
statement to the result on descending to constants in differential algebra 
\cite{bui87}.

\begin{theorem}
\label{main}
Let $S=\Spec B$ and  $X$ be a proper scheme flat over $S$ with a $q$-power
lift of Frobenius $\phi$ on $\Ou_X$ which is also perfect.
Then there exists an etale neighborhood $S' \map S$ of 
$\mfrak{p} \in S$ such that if $X' = X \times_S S'$, then $X'\simeq X_0 
\times_{\Spec{k}} S'$ where $X_0 = X\times_S \Spec{k}$ is the closed 
fiber of $X$ at $\mfrak{p}$.
\end{theorem}

{\bf Acknowledgements.} We would like to thank A. Buium for many 
in-depth discussions and also introducing this topic to the author. 
We would also like to thank J. Borger for various important comments to
this paper.

\section{Review of Witt Vectors}
Witt vectors over a general Dedekind domain with finite residue fields
 was developed in \cite{bor11a}. 
For the sake of our article, we will briefly review the general construction.
Let $B$ be a Dedekind domain and fix a maximal ideal $\mfrak{p} \in \Spec B$
with $k:=B/\mfrak{p}$ a finite field and let $q=|k|$. 
Let $R$ be the $\mfrak{p}$-adic completion of $B$. Denote by $\mfrak{m}$ the 
maximal ideal of the complete, local ring $R$ and $\iota: B \inj R$ the natural
inclusion. Then let $\pi \in B$ be such that $\iota(\pi)$ generates the maximal 
ideal $\mfrak{m}$ in $R$. Since $\iota$ is an injection, by 
abuse of notation, we will consider $\pi$ as an element of $B$ as well.
Then $k\simeq R/(\pi)$. Do note here that the identity 
map on $R$ lifts the $q$-power Frobenius on $R/(\pi)$. We will now review
the theory of $\pi$-typical Witt vectors over $R$ with maximal ideal 
$\mfrak{m}$. All the rings in this section are $R$-algebras.

Let $C$ be an $A$-algebra with structure map $u:A \map C$. In this paper,
any ring homomorphism $\psi: A \map C$ will be called the {\it lift of 
Frobenius} if it satisfies the following:

(1) The reduction mod $\pi$ of $\psi$ is the $q$-power Frobenius, that is,
$\psi(x) \equiv u(x)^q$ mod $ \pi C$.

(2) The restriction of $\psi$ to $R$ is identity.

Let $C$ be an $A$-algebra with structure map $u:A \map C$. A $\pi$-derivation
$\d$ from $A$ to $C$ means a set theoretic map satisfying the following for
all $x,y \in B$

\beqar
\label{der}
\d(x+y) &=& \d (x) + \d (y) +C_\pi(u(x),u(y)) \\
\d(xy) &=& u(x)^q \d (y) +  u(y)^q \d (x) + \pi \d (x) \d (y) \\
\eeqar
such that $\d$ when restricted to $R$ is $\d(r) = (r-r^q)/\pi$ for
all $r \in R$ and 
$$C_\pi(X,Y) = 
 \left\{\begin{array}{ll} 0, & \mb{ if $R$ is positive characteristic}\\
\frac{X^q+Y^q - (X+Y)^q}{\pi}, & \mb{ otherwise}
\end{array} \right. $$
It follows that the map $\phi: A \map C$ defined as 
$$
\phi(x) := u(x)^q + \pi \d (x)
$$ 
is an $\hA$-algebra homomorphism and is a lift of the Frobenius.
On $R$, the $\pi$-derivation $\d$ associated to $\phi$ is 
given by $\d x = \frac{\phi(x) - x^q}{\pi}$.
Considering this operator $\d$ leads 
to Buium's theory of arithmetic jet spaces \cite{bui95,bui00,
bui09}.

Note that this definition depends on the choice of uniformizer $\pi$, but in a transparent way:
if $\pi'$ is another uniformizer, then $\d(x)\pi/\pi'$ is a $\pi'$-derivation, and this correspondence
induces a bijection between $\pi$-derivations and $\pi'$-derivations.

We will present three different but equivalent point of views of Witt 
vectors:
\begin{enumerate}
\item Given an $R$-algebra $A$, the ring of $\pi$-typical Witt vectors 
$W(A)$ can be defined as
the unique $R$-algebra $W(A)$ with a $\pi$-derivtion $\d$ on 
$W(A)$ such that, given any $R$-algebra $C$ with a 
$\pi$-derivation $\d$ on it and an $R$-algebra map $f:C \map
A$, there exists an unique $R$-algebra homomorphism $g:C \map W(A)$ satisfying-
$$\xymatrix{
W(B) \ar[d] & \\
A & C \ar[l]_f \ar[ul]_g 
}$$
and $g$ satisfies $g \circ \d = \d \circ g$.
In \cite{bor11a} (following the approach of \cite{joyal} to the usual 
$p$-typical Witt vectors),
the existence of such a $W(A)$ is shown and that it is also 
obtained from the classical definition of Witt vectors using ghost vectors.

\item However, if we only restrict to flat $R$-algebras, then the Witt vectors
may be classified with the universal property of the lift of Frobenius as
follows- 
given a flat $R$-algebra $A$, the ring of $\pi$-typical Witt vectors 
$W(A)$ can be defined as
the unique flat $R$-algebra $W(A)$ with a lift of $q$-power Frobenius 
$F:W(A) \map W(A)$ on it such that, given any flat $R$-algebra $C$ with a 
lift of $q$-power Frobenius $\phi$  on it and a $R$-algebra 
map $f:C \map
A$, there exists an unique $R$-algebra homomorphism $g:C \map W(A)$ satisfying-
$$\xymatrix{
W(A) \ar[d] & \\
A & C \ar[l]_f \ar[ul]_g 
}$$
and $g$ satisfies $g \circ \phi = F \circ g$.

\item Here we will review the ghost vector definition of Witt vectors.

For any $R$-algebra $A$, define the $n+1$-fold product as 
$\Pi_n A = A \times \cdots  \times A$ 
and the infinite product $\Pi_\infty A= A \times A \times \cdots$. 
Then for all $n \geq 1$ there exists an $R$-algebra map $T_w:\Pi_n A \map 
\Pi_{n-1}A$ given by $T_w(w_0,...,w_n)= (w_0,...,w_{n-1})$.
For all $n \geq 1$ define the left shift operator $F_w:\Pi_n A \map \Pi_{n-1}
A$ as $F_w(w_0,...,w_n) = (w_1,...,w_n)$

{\it A priori} consider the following just as a product of sets 
$W_n(A):=A^{n+1}$ and the map $w:W_n(A) \map \Pi_n A$ given by 
$w(x_0,...,x_n)= (w_0,...,w_n)$ where 
\begin{align}
w_i = x_0^{q^i}+ \pi x_1^{q^{i-1}}+ \cdots + \pi^i x_i.
\end{align}
The map $w$ is known as the 
{\it ghost} map. We define the $\mfrak{p}$-typical (or $\pi$-typical) Witt
vectors $W_n(A)$ by the following theorem 
\begin{theorem}
\label{wittdef}
For each $n \geq 0$, there exists a unique ring structure on $W_n(A)$ such that
$w$ becomes a natural transformation of functors of rings.
\end{theorem}
The proof of this theorem is similar to the $p$-typical case showed 
in \cite{hessl05}.
\end{enumerate}

Now we recall some important maps of the Witt vectors. Let $W_n(A)$ denote the
truncated Witt vectors of length $n+1$. Then for every $n \geq 1$ there is a 
{\it restriction} map $T:W_n(A) \map W_{n-1}(A)$ given by $T(x_0,...,x_n) = 
(x_0,..., x_{n-1})$. Note that $T$ makes $\{W_n(A)\}$ an inverse system 
and we have $W(A) = \lim_{\leftarrow} W_n(A)$.

For every $n \geq 1$, the {\it Frobenius} ring homomorphism 
$F:W_n(A) \map W_{n-1}(A)
$ can be described in terms of the ghost vector. 
The Frobenius $F:W_n(A) \map W_{n-1}(A)$ is the unique map that makes the 
following diagram commutative in a functorial way
\begin{align}
\xymatrix{
W_n(A) \ar[r]^w \ar[d]_F & \Pi_n A \ar[d]^{F_w} \\
W_{n-1}(A) \ar[r]_-w & \Pi_{n-1} A^n
} \label{F}
\end{align}


For all $n \geq 0$, we have the multiplicative Teichmuller map $\theta:A \map 
W_n(A)$ given by $x\mapsto (x,0,0,...)$. If $A$ is of positive 
characteristic, then $\theta$ is also additive and hence $\theta$ becomes an 
$\bF_q$-algebra homomorphism.

Given an $R$-algebra $C$ with a $\pi$-derivation $\d$ on it and a $f:C \map 
A$, we will now descibe the universal map $g:C \map W(A)$. 

It is enough to show in the case when both $A$ and $C$ are flat over $R$. 
In that case the ghost map $w:W(A) \map \Pi_\infty A$ is injective.
Consider the map $[\phi]: C \map \Pi_\infty C$ given by $x \mapsto 
(x,\phi(x), \phi^2(x),...)$. Let for any $R$-algebra $D$, $F_w: 
\Pi_\infty D \map 
\Pi_\infty D$ be the left shift operator, defined by $F_w(d_0,d_1,
...) = (d_1,d_2,...)$.

$$\xymatrix{
& & C \ar[ld]_{f\circ [\phi]} \ar[d]^{[\phi]}\ar@/_1pc/[lld]_g \\
W(A) \ar[r]^w \ar[d]^F & \Pi_\infty A \ar[d]^{F_w} & \Pi_\infty C 
\ar[l]_f \ar[d]^{F_w}\\
W(A) \ar[r]^w & \Pi_\infty A & \Pi_\infty C \ar[l]_f
}$$

Then by \cite{bor11a}, the map $f\circ [\phi]:C \map \Pi_\infty A$ lifts
to $W(A)$ as our universal map $g:C \map W(A)$. It is also clear from the above
diagram that $g \circ \phi = F \circ g$. Let us now give an inductive 
description of the map $g$. Let $g(x)= (x_0,x_1,\cdots) \in W(A)$. 
Then from the 
above diagram $w\circ g = f\circ [\phi]$. Therefore for all $n \geq 0$ we have
\begin{equation}
\label{Ps}
x_0^{q^n}+ \pi x_1^{q^{n-1}}+ \cdots + \pi^n x_n = f(\phi^n(x))
\end{equation}
Note that clearly $x_0 = f(x)$ and $x_1 = f(\d(x))$. 
In the case when $A$ has a $\pi$-derivation on it, set $C=A$
and $f= \mathbbm{1}$ and let us denote the induced universal map by 
$\exp_\d:=g: A \map W(A)$ given by $\exp_\d(x)= (P_0(x),P_1(x),\cdots )$.
Composing with the restriction map $T:W(A) \map 
W_n(A)$ for each $n \geq 0$, we obtain $\exp_\d:A \map W_n(A)$. 

If $A$ is an $R$-algebra with $f:R \map A$ be the structure map, since the
identity map $\mathbbm{1}$ is a $q$-power lift of Frobenius on $R$, 
we have the universal map $\exp_\d:R \map W_n(R)$ for each $n$.
Then $W_n(A)$ is also canonically an $R$-algebra by the following composition
\begin{align}
R \stk{\exp_\d}{\longrightarrow} W_n(R) \stk{W_n(f)}{\longrightarrow} W_n(A)
\end{align}

\section{The Analogue of Taylor expansion over Witt Vectors}

We will review the case of differential algebra to motivate the analogue of 
Taylor expansion in the arithmetic case of Witt vectors. Let $X= \Spec B$ 
be an affine smooth curve over $\mathbb{C}$ with a derivation $\partial$ on it.
Let $\mfrak{p} \in X$ be a closed point on it and let $u:B \map B_0:=B/
\mfrak{p} = \mathbb{C}$ be the evaluation map at $\mfrak{p}$ and denote 
$x(\mfrak{p}):= u(x)$.

Let $D_n(B):= B[t]/(t^{n+1})$ be the ring of truncated polynomials of length
$n+1$ and $D(B):= \lim_{\leftarrow} B[t]/(t^{n+1}) = B[[t]]$. Since $B$ has
a derivation $\partial$ on it, by universal property, consider the 
Hasse-Schmidt exponential map $\exp_\d:B {\map} D(B)$ 
given by
\begin{align}
\exp_\partial(x) = \sum_{i=0}^\infty \frac{\partial^{(i)}x}{i!} t^i
\end{align} 
Then note that the Taylor expansion of any function $x \in B$ about $\mfrak{p}$
and along the derivation $\partial$ can be realised as the map 
$\overline{\exp}_\partial: B \map D(B_0)= \bb{C}[[t]]$ given by the following 
composition
\begin{align}
\label{taylor}
B \stk{\exp_\partial}{\longrightarrow} D(B) \stk{D(u)}{\longrightarrow} 
D(B_0) &   \nonumber \\
 x \longmapsto \sum_{i=0}^\infty \frac{\partial^{(i)}x}{i!} t^i 
\longmapsto  \sum_{i=0}^\infty &\frac{(\partial^{(i)}x)(\mfrak{p})}{i!} t^i
\end{align}

Recall that a derivation $\partial: B \map B$ is a ring homomorphism $B\map 
D_1(B)$, the ring of truncated polynomials of length $2$. And in the 
arithmetic case, a $\pi$-derivation $\d:B \map B$ is in fact a ring
homomorphism $B \map W_1(B)$.
The Witt vectors in this arithmetic context plays a similar role as the 
truncated polynomials in the differential algebra case. Inspired by this 
analogy, we will now introduce the arithmetic analogue of the Taylor 
expansion principle.

Let $A$ be a $B$-algebra
 which has a lift of the $q$-power Frobenius $\phi$ which when restricted to
$B$ is the identity. Also consider
the evaluation map $u:A \map A_0:= A/\mfrak{p}A$ and denote $u(a)=\overline{a}$.
Then, by analogy with
\ref{taylor} we define the {\it arithmetic Taylor expansion} of $A$ about 
$\mfrak{p}$
and with respect to $\d$ to be $\overline{\exp}_\d:A \map W(A_0)$ 
given by the following composition
\begin{align}
\label{arithexp}
A \stk{\exp_\d}{\longrightarrow} W(A) \stk{W(u)}{\longrightarrow} W(A_0)
\nonumber \\
\overline{\exp}_\d(x) = (\overline{P_0(x)},\overline{P_1(x)},...)
\end{align}
Note that $P_0(x)=x$ and $P_1(x)= \d x$.

\section{Injectivity of the Arithmetic Taylor Expansion}
Let $R$ be a complete discrete valuation ring with maximal ideal $m=(\pi)$ 
where $\pi \in R$ is a generator and also further assume that the residue
field $k=R/(\pi)$ is a finite field of order $q=p^h$ for some $h$ for a 
fixed prime $p$. Let $A$ be an $R$-algebra which is 
seperated and flat over $R$ with a $\pi$-derivation $\d$ on it. Also for any
$x \in A$, let $v_\pi(x)$ denote the $\pi$-adic valuation of $x$.

For any $R$-algebra $C$, define $C_n= C/\pi^{n+1}C$ for all $n \geq 0$.
We define the subset of constants,
denoted by, $A^\d:= \{a \in A ~|~ \d a = 0\}$. For each $n$, define the subset 
$T_n:=\{x \in W_n(A)~ |~ x=(x_0,0,...,0), x_0 \in A \} \subset W_n(A)$. 
Let $T=\lim_\leftarrow T_n$. Then in the case when the $\mb{char } A \ne 
\mb{char } k$, $T$ is only closed under multiplication but not addition. 
However when $A$ is of equal positive characteristic, $T$ is a subring of $W(A)$
since the Teichmuller map $\theta:A \map W(A)$ is then a ring homomorphism.

For $n \geq 1$, let $I_n \subset W(A)$ be the ideal defined by 
$$I_n:=\{x \in W(A_0)~|~ x=(x_0,x_1, \cdots ), x_i \in A_0 \mb{ and }
x_0= x_1= \cdots =x_{n-1} = 0\}.$$
Also let $\overline{\exp}_\d: A \map W(A_0)$ be as in \ref{arithexp}, given
by $\texp(x) = (\overline{P_0(x)},\overline{P_1(x)},...)$.

\begin{lemma}
\label{one}
Let $A$ be as above. If  $v_\pi(a) \geq 1$, then $v_\pi(\d a) = v_\pi(a)-1$. 
Moreover if $\d a = 0$ and $v_\pi(a) \geq 1$ then $a=0$.
\end{lemma}
{\it Proof.} If $a =0$, then there is nothing to show. Let $a= \pi^n b$ where
$n= v_\pi(a)$ and hence $ \pi \nmid b$ and $n \geq 1$. 
Then $\d a = \pi^{n-1} (\phi(b) - \pi^{(q-1)n}b^q)=0$ 
But $v_\pi(\phi(b))=0$ and
hence $v_\pi(\d a) = v_\pi(a) -1$.

Since $v_\pi(b) = 0$ implies $v_\pi(\phi(b))=0$ . If $\d a =0$ and $a \ne 0$ 
then since $A$ is flat over $R$, that implies $\phi(b) = \pi^{(q-1)n}b^q$ 
which is a contradiction to $a \ne 0$ and we are done.  $\qed$

\begin{lemma}
\label{min1}
Let $n \leq m$. Then the function $H(x)= q^{n-x-1}(m-x)+ x-1$ is a strictly
decreasing function in the interval $0 \leq x \leq n-1$.
\end{lemma}
{\it Proof.} Differentiating $H$ with respect to $x$ we get: 
$$H'(x) = -q^{n-1-x}[1+ (\ln{q})(m-x)] +1$$ 
Note that $1+ (\ln{q})(m-x) >1$ for all $x \leq m-1$ and hence 
$-q^{n-1+x}[1+ (\ln{q})(m-x) < -1$ which implies $H'(x) < 0$ 
and we are done. $\qed$

\begin{proposition}
\label{explicit}
For all $n$ we have 
\begin{align}
P_n(x)= \sum_{i=0}^{n-1}\sum_{j=1}^{q^{n-1-i}} \pi^{i+j-n} 
\left(\begin{array}{l} q^{n-1-i}\\ j\end{array}\right)
P_i(x)^{q(q^{n-1-i}-j)}(\d P_i(x))^j
\end{align}
\end{proposition}
{\it Proof.} 
From \ref{Ps}, we get-
\begin{equation}
x^{q^n} + \pi P_1(x)^{q^{n-1}}+ \pi^2 P_2(x)^{q^{n-2}}+ ... + \pi^nP_n(x) =
\phi^n(x)
\end{equation}
Since $A$ has  a $\pi$-derivation we have
\beqar
\pi^nP_n(x) + \sum_{i=0}^{n-1} \pi^i P_i(x)^{q^{n-i}} &=& \phi(\phi^{n-1}(x)) \\
&=& \phi(\sum_{i=0}^{n-1} \pi^i P_i(x)^{q^{n-1-i}})\\
&=& \sum_{i=0}^{n-1} \pi^i \phi (P_i(x))^{q^{n-1-i}} \\
P_n(x) &=& \sum_{i=0}^{n-1} \pi^{i-n} (\phi(P_i(x))^{q^{n-1-i}} - 
P_i(x)^{q^{n-i}}) \\
\eeqar
and the result follows by binomially expanding 
$\phi(P_i(x))^{q^{n-1-i}}= (P_i(x)^q+ \pi \d P_i(x))^{q^{n-1-i}}$, for each 
$i$. $\qed$

\subsection{The unequal characteristic case}
Let $q=p^h=|k|$ as before and $v_\pi(p)= e$ be the absolute 
ramification index of $p$ in $R$, that is $pR= (\pi)^e$.
Let us define the following
\begin{align}
L_{ij} &= \pi^{i+j-n} \left(\begin{array}{l} q^{n-1-i}\\ j\end{array}\right)
P_i(x)^{q(q^{n-1-i}-j)}(\d P_i(x))^j \nonumber\\
S_i &= \sum_{j=1}^{q^{n-1-i}} L_{ij}
\end{align}
Then 
\begin{align}
P_n(x) = \sum_{i=0}^{n-1} S_i
\end{align}

\begin{lemma}
\label{l1}
If $v_\pi(P_i(x)) = m-i,~ 0 \leq i < m$, then $v_\pi(L_{ij}) = 
i-n +(n-1-i)eh + (m-i)q^{n-i} - (m-i)(q-1)j-v_\pi(j)$.
\end{lemma}
{\it Proof.} We know that $v_\pi\left( \left(\begin{array}{l} p^l\\ j \end{array}
\right) \right) = lv_\pi(p) - v_\pi(j)$. Then the result follows easily from the
following computation
\beqar
v_\pi(L_{ij}) &=& i+j-n + (n-1-i)eh - v_\pi(j) + (m-i)q(q^{n-1-i}-j) +
(m-i-1)j  \\
&=& i-n +(n-1-i)eh + (m-i)q^{n-i} - (m-i)(q-1)j -v_\pi(j) ~~ \qed\\
\eeqar

\begin{lemma}
\label{l2}
If $v_\pi(P_i(x)) = m-i,~0\leq i < m$, then $v_\pi(S_i)= v_\pi(L_{iq^{n-1-i}})
=q^{n-1-i}(m-i)-n+i$.
\end{lemma}
{\it Proof.} Since $A$ is flat over $R$, it is sufficient to show that 
$v_\pi(L_{iq^{n-1-i}}) \lneq v_\pi(L_{ij}$ for all $1 \leq j \leq q^{n-1-i}-1$.
Let $C = (m-i)(q-1) > 0$. Then for all $j,j'= 0,...,q^{n-1-i}$ by lemma 
\ref{l1} we have 
$v_\pi(L_{ij'}) -v_\pi(L_{ij})= C(j-j') +(v_\pi(j)-v_\pi(j'))$

Now for all $0 \leq j \le q^{n-1-i}-1$ we have 
\begin{align}
C(q^{n-1-i}-j) > 0 > v_\pi(j) - v_\pi(q^{n-1-i}) \nonumber\\
0 > C(j-q^{n-1-i}) + (v_\pi(j)-v_\pi(q^{n-1-i})) \nonumber \\
 = v_\pi(L_{iq^{n-1-i}}) - v_\pi(L_{ij})
\end{align}
and we are done. $\qed$

Now we return to the general characteristic case
\begin{theorem}
\label{val}
Let $A$ be a separated flat $R$-algebra with a $\pi$-derivation $\d$ on it
and $x \in A$ with $v_\pi(x)=m < \infty$ (in particular $x \ne 0$).
Then for all $n \leq m$, $P_n(x) \ne 0$ and moreover 
$$v_\pi(P_n(x)) = m-n.$$
\end{theorem}
{\it Proof.} We will prove the result by induction on $n$. For $n=1$, we
have already shown that in lemma \ref{one}. Assume it is true for n-1. Then
$v_\pi(P_i(x)) = m-i$ for all $i=0,...,n-1$.

We wish to first determine the valuation $v_\pi(S_i)$ for all $i=0,...,n-1$.
Let us consider two separate cases of $\mb{char }A =p$ and $\mb{char }A
\ne p$.
If $\mb{char }A=p$, then note that $S_i = L_{iq^{n-1-i}}$ and therefore 
$v_\pi(S_i)=v_\pi(L_{iq^{n-1-i}})$. Otherwise by lemma \ref{l2} we also 
have $v_\pi(S_i) = v_\pi(L_{iq^{n-1-i}})$.

Therefore for all $i=0,..., n-1$, we have
\beqar
v_\pi(S_i) &=& q^{n-1-i}(m-i) -(n- i)
\eeqar
By lemma \ref{min1} the right hand side is a strictly decreasing 
function of $i$ and hence the minimum is attained at $i=n-1$. Since
$A$ is flat, we conclude that 
\beqar
v_\pi(P_n(x)) &=& \min_{i=0,...,n-1} {(q^{n-1-i}(m-i)+ i-n)}\\
&=& m-n
\eeqar
and we are done. $\qed$

{\it Proof of Theorem \ref{nokernel}.} 
Suppose $x \in \ker(\overline{\exp}_\d)$ and $x \ne 0$. 
Then $\overline{P_i(x)} = 0$ for all $i \geq 0$. 
Since $\overline{P_0(x)} =\overline{x}$, clearly $v_\pi(x) \geq 1$ and let 
$n:=v_\pi(x) \geq 1$. However, by theorem \ref{val}, 
$v_\pi(P_n(x)) =0$ which implies $\overline{P_n(x)} \ne 0$ and hence is a
contradiction and therefore we must have $x=0$ and we are done.$\qed$

\begin{corollary}
\label{modinj}
The induced map $A_n \map W_n(A_0)$ is an injection.
\end{corollary}
{\it Proof.} Consider the truncated map $\texp:A \map W_n(A_0)$ given by 
\newline
$x \mapsto (\overline{P_0(x)},..., \overline{P_n(x)})$. Let $x$ 
be in the kernel of the map. Then by theorem \ref{val}, $v_\pi(x) \geq n+1$,
that is $x \in \pi^{n+1}A$. And clearly $\pi^{n+1}A$ is contained inside the
kernel by the same theorem and hence we are done. $\qed$

Consider the topology on $W(A_0)$ induced by the basis of open sets $I_n$ 
for all $n$. Then the above corollary implies the following about the $\pi$-adic
topology of $A$:

\begin{corollary}
\label{topology}
Let $A$ be as above. The $\pi$-adic topology on $A$ is the restriction of the
$I_n$-adic topology on $W(A_0)$ by the map $\overline{\exp}_\d$.
\end{corollary}

\begin{lemma}
\label{same}
Let $C \inj D$ be two $\pi$-adically complete, flat, $R$-algebras
such that $C_0 \simeq D_0$. Then $C\simeq D$. 
\end{lemma}
{\it Proof.} We prove $C_n \simeq D_n$ using induction. 
Since both $C$ and $D$ are $\pi$-torsion free, for all $n$, we have 
$\pi^{n-1}C/\pi^nC \simeq C_0 \simeq D_0 \simeq \pi^{n-1}D/\pi^nD$.
For $n=0$, it is the
hypothesis. Suppose true for $n-1$. Then
$$\xymatrix{
0 \ar[r]& \pi^{n-1}C/\pi^n C \ar[d] \ar[r] & C_n \ar[d] \ar[r] & C_{n-1} \ar[d]
\ar[r] & 0 \\
0 \ar[r]&  \pi^{n-1}D/\pi^n D \ar[r] & D_n \ar[r] & D_{n-1} \ar[r] & 0 \\
}$$
and we conclude by snake lemma. $\qed$

The following lemma is standard.
\begin{lemma}
\label{Bhat}
Let $B,\pi$ and $k$ be as above. Also further assume $B$ is a $k$-algebra. Then
$R \simeq k[[\pi]]$.
\end{lemma}
{\it Proof.} Since $B$ is a $k$-algebra implies $R$ is a $k$-algebra.
We claim that the image of $k[[\pi]]$ in $R$ is injective. 
Let $f \in k[[\pi]]$ be such that its image in $R$ is $0$. Then there exists
an $l \geq 0$ such that $f$ can be written as
$f= \sum_{i=0}^\infty \alpha_{l+i} \pi^{l+i}$, where $\alpha_{l+i} \in k$ for 
all $i$ and $a_l \ne 0$.
Since $R$ is $\pi$-torsion free, that implies $\alpha_l = \pi(\alpha_{l+1}+
\cdots)$ which implies $v_\pi(\alpha_l) \geq 1$. But since $\alpha_l \in k$ we 
must have $\alpha_l=0$ which is a contradiction and this proves the claim.
Now we have $k[[\pi]] \subseteq R$ and both are $\pi$-adically complete and
has the same residue field $k$. This implies they are isomorphic by lemma 
\ref{same} and we are done. $\qed$

If $A$ is of equal positive characteristic, then since the Witt vector addition
is component wise linear, it makes $T$ a subgroup under addition too and hence 
is a subring of $W(A_0)$.  Therefore $A^\d$, which is isomorphic to $T$,
is also a  subring of $A$. 

\begin{lemma}
\label{allzero}
If $\d x =0 $ then $P_n(x)= 0$ for all $n \geq 1$.
\end{lemma}
{\it Proof.} 
We will prove by induction on $n$. It is true for $n=1$ by hypothesis. If 
true for $n-1$ then by proposition \ref{explicit}, $P_n(x)$ can be expressed as
$$P_n(x)= \sum_{i=0}^{n-1}\sum_{j=1}^{q^{n-1-i}} \pi^{i+j-n} 
\left(\begin{array}{l} q^{n-1-i}\\ j\end{array}\right)
P_i(x)^{q(q^{n-1-i}-j)}(\d P_i(x))^j
$$
which clearly implies $P_n(x)=0$ and we are done. $\qed$

\begin{lemma}
\label{wo}
For all $1 \leq x \leq n-1 \leq m$, $ H(x)= x -n + q^{n-1-x}(m-x+1)$ is a
strictly decreasing function in $x$.
\end{lemma}
{\it Proof.} Differentiating $H$ we get 
$$H'(x) = 1 - q^{n-1-x}((m-x+1) \ln q +1)$$
Note that $q^{n-1-x} \geq 1$ and $(m-x+1)\ln q+1 > 1$ for all $1 \leq x 
\leq n-1 \leq m$ and therefore $H'(x) < 0$ and we are done. $\qed$

\begin{proposition}
\label{modpiP}
Let $R$ be of equal characteristic and $x \in A$ be such that $\vpi (x) = 0$ 
and $\vpi(\d x) =
m \geq 1$. Then $\vpi(P_n(x)) = m -n +1$ for all $1 \leq n \leq m+1$.
\end{proposition}
{\it Proof.} We have  for all $n$,
$$P_n(x) = \sum_{i=0}^{n-1} S_i$$ 
where $S_i = \pi^{i-n}(\phi(P_i(x))^{q^{n-1-i}} - P_i(x)^{q^{n-i}})$ for all
$i=0 \cdots n-1$.
We will prove our result by induction. For $n=1$, $P_1(x)= \d x$ and hence
the result is true by hypothesis. Now assume true for $n-1$.
Since $\ch R >0 $ we have $S_0 =\phi(x)^{q^{n-1}} -x^{q^n}= (x^{q^n} + 
(\pi\d x)^{q^{n-1}}) -x^{q^n} =
(\pi \d x)^{q^{n-1}}$ which implies, $\vpi (S_0) = q^{n-1}(1+m)$.

Now for all $1 \leq i \leq n-1 \leq m$, we have $\vpi (\phi(P_i(x)^{q^{n-1-i}})
 < \vpi ( P_i(x)^{q^{n-i}})$ and therefore we have
$$ \vpi(S_i)= i-n + q^{n-1-i}(m-i+1) $$
Then by lemma \ref{wo}, the valuation $\vpi(S_i)$ is a strictly decreasing
function of $i$ and therefore we have 
$$\vpi\left(\sum_{i=1}^{n-1} S_i\right) = \vpi(S_{n-1}) = m-n+1$$
Note that $\vpi(S_0)= q^{n-1}(1+m) > m-n+1 = \vpi(\sum_{i=1}^{n-1} S_i)$ for 
all $n > 1$. Therefore we have 
$$\vpi(P_n(x)) = \vpi \left(\sum_{i=0}^{n-1} S_i\right) = 
\vpi\left(\sum_{i=1}^{n-1} S_i \right) = m-n+1$$
and we are done. $\qed$

\begin{proposition}
\label{sadhu}
If $R$ is of positive characteristic, then $\oexp(A^\d) = \oexp(A) \cap T$.
\end{proposition}
{\it Proof.}
Let $x \in A^\d$. Then since $\d x = 0$, by lemma \ref{allzero}, we have
$P_n(x) = 0$ for all $n \geq 1$ and hence $\oexp(A^\d)
 \subseteq \oexp(A) \cap T$.

Conversely, consider an element $y \in \oexp(A) \cap T$. Then there exists an 
$x \in A$ such that $\oexp(x) = y$ and satisfies $\overline{P_n(x)} = 0$
for all $n \geq 1$. Since $\oexp$ is injective, it is enough to consider the
case when $\overline{x} \ne 0$. Now if $\d x =0$ then we are done. Otherwise
if $\d x \ne 0$  then we have $0< \vpi(\d x) =m < \infty$, since 
$\overline{P_1(x)}= \overline{\d x} = 0$ and $A$ is separated. But by 
proposition \ref{modpiP} we have $\vpi(P_{m+1}(x)) = 0$ and hence we have
$\overline{P_{m+1}(x)} \ne 0$ which is a contradiction and therefore we must 
have $\d x = 0$ and therefore $\oexp(A) \cap T \subseteq \oexp(A^\d)$ and 
we are done. $\qed$

\begin{lemma}
If $R$ is of positive characteristic, then $A^\d \cap \pi A= (0)$. In 
particular the composition 
\begin{align}
A^\d \inj A \stk{u}{\map} A_0 \nonumber
\end{align}
is an injection. 
\end{lemma}
{\it Proof.} This follows because $\oexp(A^\d) \subseteq T$ and by proposition
\ref{sadhu}, $T \cap ~\oexp(\pi A) = (0)$ as well as the map $\oexp$ is an 
injection. $\qed$

\begin{lemma}
Let $R$ be of positive characteristic and $A$ be a separated flat $R$-algebra 
Then the subring generated by finite sums of $A^\d$ and $\pi$ inside $A$ is 
$A^\d[\pi]$.
\end{lemma}
{\it Proof.} If not there exists $0 \ne f = a_l \pi^l+ ... a_n\pi^n$ where
$a_l,...,a_n \in A^\d$ and $a_l \ne 0$ such that $f = 0$. Since $A$ is flat 
over $R$, we have $a_l= -a_{l+1}\pi+...+ a_n\pi^{n-l}\in \pi A$, that is,
 $v(a_l) \geq 1$. But since $a_l \in A^\d$ we also have $\d a_l =0$.
But that implies, by lemma \ref{one}, $a_l =0$ which is a contradiction.
Therefore we must have $f=0$ and we are done. $\qed$

\begin{corollary}
Let $R$ be of equal characteristic and $A$ be a separated flat $\pi$-adically
complete $R$-algebra.Then
the $\pi$-adic completion of the subring generated by $A^\d$ and $\pi$ inside
$A$ is $A^\d[[\pi]]$.
\end{corollary}

{\it Proof of Theorem \ref{iso}.}  
Since $\phi$ induces an isomorphism on $A_0$, we have that $\oexp:
A \simeq W(A_0)$ is an isomorphism. Therefore by proposition \ref{sadhu} 
we have $A^\d \simeq T  \simeq A_0$. Hence $A^\d[[\pi]] \inj A$ and the 
result follows from lemma \ref{same}.	$\qed$

\begin{corollary}
\label{iso-cor}
Let $A$ be a $\pi$-torsion free separated and $\pi$-adicially complete 
$R$-algebra with a lift of Frobenius $\phi$ that induces an isomorphism on 
$A_0$. Also further assume $B$ is a $k$-algebra. Then $A_n \simeq
A^\d \otimes_k B_n$.
\end{corollary}
{\it Proof.} From theorem \ref{iso} we have $A_n \simeq A^\d \otimes_k 
k[\pi]/(\pi^{n+1})$. By lemma \ref{Bhat}, we have $B_n \simeq 
k[\pi]/(\pi^{n+1})$ and we are done. $\qed$

\section{An Application in the case of Schemes}

%
%

Let $(B,\mfrak{p})$ be as before. Furthermore assume $B$ is also a $k$-algebra 
where $k= B/\mfrak{p}$. And as before, let $R$ be the completion of $B$
with respect to $\mfrak{p}$ and $\iota: B\map R$ be the canonical injective 
map. Let $A$ be a $B$-algebra and let $\hA$ denote the completion of $A$
with respect to  $\mfrak{p}A$. Then $\hA$ is an $R$-algebra. 

\begin{lemma}
\label{nopi4}
Let $A$ be $\pi$-torsion free, then $\hA$ is $\pi$-torsion free and 
separated.
\end{lemma}
{\it Proof.} Since $A$ has no $\pi$-torsion, $\hA$ is $\pi$-torsion free as
well. And $\hA$ is separated since it is the $\pi$-adic 
completion of $A$. $\qed$

{\it Proof of Theorem \ref{main}.}
Let $S_n = \Spec~ B/\mfrak{p}^{n+1}$. For any scheme $Z$ over $S$, let 
$Z_n := Z \times_S S_n$.
Consider $Y = X_0 \times_{\Spec~k} S$. 
Since $X$ is an integral scheme of finite type over $S$, 
by lemma \ref{nopi4}, we
have that $\Ou_{\hat{X}}$ is a sheaf of $\pi$-torsion free separated 
Noetherian $R$-algebras. Therefore by corollary \ref{iso-cor}, we have
for each $n$, compatible isomorphims $f_n:X_n \map Y_n$. Then by Artin 
approximation, \cite{artin69} cor 2.4, there exists an etale neighbourhood
$S'$ of $\mfrak{p} \in S$ such that if $X'=X\times_S S'$ then $X' \simeq
Y\times_S S' \simeq X_0 \times_{\Spec~k} S'$.
$\qed$ 

 

\footnotesize{

}

\end{document}